\documentclass[11pt]{article}

\usepackage{amssymb}

\usepackage{latexsym}
\newtheorem{lemma}{Lemma}[section]
\newtheorem{theorem}[lemma]{Theorem}
\newtheorem{cor}[lemma]{Corollary}
\newtheorem{prop}[lemma]{Proposition}
\newtheorem{rem}[lemma]{Remark}

\newtheorem{notn}[lemma]{Notation}
\def\ZZ{\mbox{{\rm Z\hskip-4pt \rm Z}}}   
\def\CC{{\rm \kern.24em \vrule width .02em height1.4exdepth-.05ex
\kern -.26em C}}
\newcommand{\RR}{\mbox{$I\!\!R$}}

\newfont{\tenrm}{cmr10 scaled 1000}
\title{The quotient Unimodular Vector group is nilpotent} 

\date{}
\author{{\large Reema Khanna \& ~Selby Jose \&}\\
{\large Sampat Sharma \& Ravi A. ~Rao}\\
\small{K.J. Somaiya College, Vidyavihar, Mumbai 400 077 \&}\\ 
\small{Department of Mathematics, Institute of Science, Mumbai 400 032 \&}\\
\small{School of Mathematics, Tata Institute of Fundamental Research,}\\
\small{1, Dr. Homi Bhabha Road, Mumbai 400 005.}
}
\date{}

\begin{document}

\maketitle
\begin{abstract}
Jose--Rao introduced and studied the Special Unimodular 
Vector group $SUm_r(R)$ and $EUm_r(R)$, its Elementary Unimodular Vector 
subgroup. They proved that for $r \geq 2$, $EUm_r(R)$ is a normal subgroup of $SUm_r(R)$. 
The Jose--Rao theorem says that the quotient Unimodular Vector 
group, $SUm_r(R)/EUm_r(R)$, for $r \geq 2$, is a subgroup of the orthogonal 
quotient group $SO_{2(r+1)}(R)/EO_{2(r + 1)}(R)$. The latter group is known to be 
nilpotent by the work of Hazrat--Vavilov, following methods of A. Bak; and 
so is the former. 

In this article we give a direct proof, following ideas of A. Bak, to show that 
the quotient Unimodular Vector group is nilpotent of class $\leq d 
= \dim(R)$. We also use the Quillen--Suslin theory, inspired by A. Bak's 
method, to prove that if $R = A[X]$, with 
$A$ a local ring, then the quotient Unimodular Vector group is abelian\footnote{ \S 4 is part of the doctoral thesis of the
first named author under the second named author; \S 5 is part of the doctoral thesis of the third named author under the fourth named author.}.
\end{abstract}

\section{Introduction}
$R$ will be a commutative ring with $1$, in which $2$ is invertible.  
$Um_{r + 1}(R)$ will denote the set of unimodular vectors 
$v \in R^{r + 1}$, i.e. those vectors $v$ for which there is a vector 
$w \in R^{r + 1}$, with $\langle v, w \rangle = v\cdot w^T = 1$.

Suslin introduced the Suslin matrix in (\cite{8}, \S 5), and indicated its 
properties as well as how he felt they will be useful. 
 
In \cite{4} we initiated the study of the special unimodular vector
group $SUm_r(R)$, which is a subgroup of  
$GL_{2^r}(R)$ related to $Um_{r + 1}(R)$. We also introduced the 
elementary unimodular vector subgroup $EUm_r(R)$ of $SUm_r(R)$, which is 
related to the $(r + 1)$-unimodular vectors which have a completion to an 
elementary
matrix. We developed the calculus for $EUm_r(R)$ in \cite{4}, and got a 
nice set of generators for it. In \cite{5} we showed that $EUm_r(R)$ is a 
normal subgroup of $SUm_r(R)$, for $r \geq 2$. 

In \cite{S} Suslin, inspired by Quillen's methods in \cite{Q}, applied them 
to the study of unstable $K_1$-theory of polynomial rings. He proved the 
$K_1$-analogue of the Local-Global Principle and the Monic Inversion Principle. 
The theory built up in \cite{Q,S} is known as the Quillen--Suslin theory.

Using Quillen--Suslin Local Global principle, A. Bak established in 
\cite{1}, that the linear quotient $SL_n(R)/E_n(R)$, for $n \geq 3$, is 
nilpotent. This 
theme has been revisited several times for different classical groups, 
see \cite{hv}, \cite{astep}, and (\cite{rb}, \S 3.3) for instance. 

Now we apply Bak's approach to the pair $(SUm_r(R), EUm_r(R))$, 
for $r \geq 2$, when $R$ is a noetherian ring of Krull dimension $d$. 
We give a direct approach to reprove the result in \cite{6} that the 
unimodular vector quotient $SUm_r(R)/EUm_r(R)$ is a nilpotent group of 
class $d$. (The latter had been established in \cite{6} via the Jose--Rao
theorem that the unimodular vector quotient group was a subgroup of the 
special orthogonal quotient group; which was nilpotent in view of 
\cite{hv}.) 

We also deduce a relative version of this result from the absolute case. 
This argument does not depend on the Excision ring argument of W. van der 
Kallen, which is normally used to deduce `relative' results; and is much 
more flexible. (This approach evolved from the work \cite{hom} according to 
Anjan Gupta; who used it in his thesis (\cite{ag}, \S 2.2) to reprove a 
theorem of Chattopadhyay--Rao in \cite{pcrr}).  

Finally, we consider $SUm_r(R)/EUm_r(R)$, the unimodular vector quotient group, 
when 
$R = A[X]$ is a polynomial extension of a local ring $A$. In this case 
we show, arguing as in \cite{rrss} that the unimodular quotient group is 
an abelian group. A relative version for extended ideals is also deduced.

\section{Recap about the Suslin matrix $S_r(v, w)$}

Given two row vectors $v, w \in R^{r +1}$,  
A. Suslin constructed in [\cite{8}, \S 5], a matrix 
$S_r(v, w)$, which is of determinant one if 
$\langle v, w \rangle =  v\cdot w^T = 1$. He defined this inductively, 
  as follows: Let 
$v = (a_0, a_1, \ldots, a_r) = (a_0, v_1)$,
with $v_1 = (a_1, \ldots, a_r)$, 
$w = (b_0, b_1, \ldots, b_r) = (b_0, w_1)$,
with $w_1 = (b_1, \cdots, b_r)$. Set $S_0(v, w) = a_0$, and set 
\begin{eqnarray*}
S_r(v, w)  & = &  \pmatrix{a_0 I_{2^{r - 1}} & S_{r - 1}(v_1,
w_1)\cr - S_{r - 1}(w_1, v_1)^T & b_0 I_{2^{r - 1}}\cr}.
\end{eqnarray*}

The reader will find more details about these matrices in this amazing 
$\S 5$; with several unresolved questions. 

\vskip0.1in
These matrices have been studied by Jose--Rao in \cite{5,6}. The survey 
article \cite{surv} gives a quick glimpse at the known results today. 

We shall denote by 
$SUm_r(R)$ the subgroup of $GL_{2^r}(R)$ generated by the set
$\{S_r(v, w)| v, w\in 
R^{r + 1}$, $\langle v, w \rangle = 1\}$, and $EUm_r(R)$ its subgroup 
generated by the set $\{S_r(v, w)| v, w\in 
R^{r + 1}$, $\langle v, w \rangle = 1$, $v = e_1\varepsilon$, for some 
$\varepsilon \in E_{r + 1}(R)\}$. It was shown in \cite{5}, that 
$EUm_r(R)$ is a normal subgroup of $SUm_r(R)$, for $r \geq 2$. 

For a matrix $\alpha \in M_k(R)$, we define $\alpha^{top}$
as the matrix whose entries are the same as that of $\alpha$ above the
diagonal, and on the diagonal, and is zero below the diagonal. Similarly,
we define $\alpha^{bot}$. Moreover, we use $\alpha^{tb}$ for $\alpha^{top}$ or $\alpha^{bot}.$

In \cite{4}  a structure theorem for $EUm_r(R)$ was proved. The 
following  nice set of generators of $EUm_r(R)$ was established:

For $2 \leq i \leq r+1$, $\lambda \in R$, let
$$\begin{array}{cccccc}
E(e_i)(\lambda) & = & S_r(e_1 + \lambda e_{i}, e_1), & E(e_i^*)(\lambda) & = & S_r(e_1, e_1 + \lambda e_{i}), \\
E(e_{i1})(\lambda) & = & S_r(e_i + \lambda e_1,e_i), & E(e^*_{i1})(\lambda)&=& S_r(e_i, e_i + \lambda e_1).
\end{array}$$
It was shown that the group $EUm_r(R)$ can be generated by either
\begin{enumerate}
\item[(a)] $E(c)(x)$, $E(d)(x)S_r(e_i, e_i)^{-1}$, if $2$ is invertible 
in $R$, or by 
\item[(b)] $E(c)(x)^{top}$, $E(c)(x)^{bot}$,  
\end{enumerate}
where $c = e_i$ or $e_i^*$, $d = e_{i1}$ or $e_{i1}^*$, $2 \leq i \leq 
r + 1$, $x \in R$.

In \cite{4, 6} Jose--Rao noted a fundamental property which is satisfied by the 
Suslin matrices. Let
$v$, $w$, $s$, $t \in M_{1, r+1} (R)$. Then
\begin{eqnarray*}
S_r (s, t) S_r (v,w) S_r (s,t) &=& S_r (v',w')       \\
S_r (t, s) S_r (w,v) S_r (t,s) &=& S_r (w', v'),   
\end{eqnarray*}
for some $v',w' \in M_{1r+1}(R)$, which depend linearly on $v$, $w$ and 
quadratically on $s$, $t$.  Consequently, 
$v' \cdot w^{'T}= (s \cdot t^T)^2 (v \cdot w^T)$.

This fundamental property enables one to define an involution $\star$
on the group $SUm_r(R)$, details of which can be found in \cite{6}. 
This involution 
is then used to give an action of $SUm_r(R)$ on the Suslin space, 
{\it viz.} the free $R$-module of rank $2(r + 1)$ 
$$S = \{S_r(v, w) | v, w \in M_{1 r+ 1}(R)\}.$$
(For a basis one can take $se_1, \ldots, 
se_{r + 1}, se_1^*, \ldots, se_{r + 1}^*$, where 
$se_i = S_r(e_i, 0)$,  $se_i^* = S_r(0, e_i)$, for $1 \leq i \leq r$.) 

In \cite{6} they associated a linear transformation $T_g$ of the 
Suslin space with a Suslin matrix $g$, via 
$$T_g(x,y) = (x',y'),$$ 
where $gS_r(x, y)g^* = S_r(x', y')$. Moreover, if $g$ is a product of 
Suslin matrices $S_r(v_i, w_i)$, with $\langle v_i, w_i\rangle = 1$, for all 
$i$, then $T_g \in SO_{2(r + 1)}(R)$, i.e. 
$$\langle T_g(v, w), T_g(s, t)\rangle = \langle (v,w), (s, t) \rangle =
v\cdot w^T + s\cdot t^T.$$

\section{Computation of the matrix of 
the linear transformation}

In (\cite{6}, \S 4), via the fundamental property, Jose-Rao observed that  
the above action induces a canonical homomorphism 
\begin{eqnarray*} 
\varphi: {SUm}_r(R) &\rightarrow & {SO}_{2(r + 1)}(R),\\
\varphi(S_r(v, w)) &=& T_{S_r(v, w)} =
\tau_{(v, w)} \circ \tau_{(e_1 e_1)},
\end{eqnarray*}
where $\tau_{(v, w)}$ is the standard reflection with respect to the vector 
$(v, w) \in R^{2(r + 1)}$ (of length one) given by the formula
\begin{eqnarray*}
\tau_{(v, w)}(s, t) & = & \langle v, w \rangle (s, t) - (\langle v, t \rangle
+ \langle s, w \rangle)(v, w).
\end{eqnarray*}

The following simple computation gives an alternate way to prove this:

\begin{lemma}\label{3.1}
Let $R$  be a commutative ring  with $1$. Let $v,  w \in {\rm Um}_{r+1}(R)$,
then  the  matrix of  the  linear  transformation $T_{S_r(v,w)}$  with
respect to  the (ordered) basis  $$\{S_r(e_1, 0), S_r(e_2,  0), \cdots,
S_r(e_{r+1}, 0), S_r(0,e_1),  S_r(0,e_2), \cdots, S_r(0,e_{r+1})\}$$ is
$$\left(I -  \pmatrix{v^T \cr w^T}  \pmatrix{w & v} \right)  \left(I -
\pmatrix{e_1^T   \cr  e_1^T}   \pmatrix{e_1  &   e_1}   \right).$$  
\end{lemma}

\noindent
Proof: Let  $v =  (a_0, a_1,  \cdots, a_r)$, $w  = (b_0,  b_1, \cdots,
b_r)$.  By the definition of $T_{S_r(v,w)}$,
\begin{eqnarray*}
T_{S_r(v,w)}(e_1, 0) & =  & \tau_{(v,w)} \circ \tau_{(e_1, e_1)} (e_1,
0)\\ & = & \tau_{(v,w)}(0, -e_1) = (0, -e_1) + a_0 (v,w) = (a_0v, a_0w
- e_1)\\ 
 T_{S_r(v,w)}(e_j, 0)  & = & \tau_{(v,w)} \circ \tau_{(e_1,
e_1)} (e_j, 0)\\ & = & \tau_{(v,w)}(e_j,0) = (e_j,0) - b_{j-1} (v,w) =
(e_j  -   b_{j-1}v,  -b_{j-1}w)\\  
 T_{S_r(v,w)}(0,  e_1)   &  =  &
\tau_{(v,w)}    \circ   \tau_{(e_1,    e_1)}   (0,e_1)\\    &    =   &
\tau_{(v,w)}(-e_1,0) = (-e_1,0) + b_0  (v,w) = (b_0v - e_1, b_0w)\\ 
T_{S_r(v,w)}(0,e_j)  &  =   &  \tau_{(v,w)}  \circ  \tau_{(e_1,  e_1)}
(0,e_j)\\ &  = &  \tau_{(v,w)}(0,e_j) = (0,e_j)  - a_{j-1} (v,w)  = (-
a_{j-1}v, e_j - a_{j-1}w)
\end{eqnarray*}
Thus the matrix of $T_{S_r(v,w)}$ is
$${\small \pmatrix{a_0v & e_2 - b_1v & \cdots  & e_{r+1} - b_rv & b_0v - e_1 &
-a_1v & \cdots & -a_rv\cr a_0w - e_1 & -b_1w & \cdots & -b_rw & b_0w &
e_2 -  a_1w &  \cdots &  e_{r+1} - a_rw}.}$$  

Right multiply  the above
matrix by the matrix $\left(I - \pmatrix{e_1^T \cr e_1^T} \pmatrix{e_1
& e_1} \right)$ will interchange  the 1-st and $(r+2)$-th columns with
sign changed.  Hence, the matrix of $T_{S_r(v,w)}$ is
 $$\left(I - \pmatrix{v^T  \cr w^T} \pmatrix{w & v}  \right) \left(I -
\pmatrix{e_1^T \cr  e_1^T} \pmatrix{e_1 & e_1}  \right)$$ as required.
\hfill$\Box$

\vskip0.15in

\noindent
{\bf Notation:} We denote the matrix of the linear transformation 
$T_{S_r(v,w)}$ by $[T_{S_r(v,w)}]$. 

\vskip0.1in

Let us recollect the matrix of the linear transformations corresponding to 
the generators of $EUm_r(R)$, $r\geq 2$, computed in \cite{6}.

\medskip

For the sake of completeness we give a slightly simpler argument than the 
one given in \cite{6} below. However, in this approach, unlike in 
\cite{6},  we need that $2$ is invertible in $R$. 

\begin{lemma}
For $2 \le i, j \le r+1$, one has the following relations in EUm$_{r}(R)$: 
\begin{eqnarray*}
E(e_i^*)(-2\lambda)^{bot}   &  =  &   S_r(e_1  -   e_j,  e_1   -  e_i)
S_r((1+\lambda)e_1 + e_j, e_1 - \lambda e_j) \\ 
& & S_r(e_1 - e_j, e_1+   e_i)S_r((1-\lambda)e_1 + e_j, e_1 + \lambda   e_j)\\
& & \left[E(e_j^*)(\lambda), E(e_i^*)(1)\right].\\
E(e_i)(-2\lambda)^{bot} & = & \left[E(e_i)(-1), E(e_j)(-\lambda)\right]\\
&&S_r(e_1+\lambda e_j, (1 - \lambda)e_1 + e_j) S_r(e_1+e_i, e_1+e_j)\\
&&S_r(e_1-\lambda e_j, (1+\lambda)e_1 + e_j) S_r(e_1-e_i, e_1-e_j).
\end{eqnarray*}
$($Note that by reversing the elements in the 
product in the above relation we can obtain the formulae for 
$E(e_i^*)(-2\lambda)^{top}$ and  $E(e_i)(-2\lambda)^{top}$.$)$
\end{lemma}

\noindent
Proof: 
We prove the first relation; the others are verified similarly. 
Put  $x =1$, $y  = \lambda$, and  $z = 1$  in the proof  of [\cite{4}, Proposition 5.6], to get
\vskip.2in
\noindent
$E(e_i^*)(-2\lambda)^{bot}  $
\begin{eqnarray*}
&   =  &   \{E(e_j)(1)^{-1}\}\{E(e_j)(1/2)
E(e_i^*)(1/2)^{-1}    E(e_i^*)(1/2)^{-1}E(e_j)(1/2)   \}   \\    &   &
\{E(e_j)(1)^{-1}\}  \{S_r((1+\lambda)e_1 + e_j,  e_1 -  \lambda e_j)\}
\{E(e_j)(1)^{-1}\}\\            &            &           \{E(e_j)(1/2)
E(e_i^*)(1/2)E(e_i^*)(1/2)E(e_j)(1/2)  \}   \{E(e_j)(1)^{-1}\}\\  &  &
\{S_r((1-\lambda)e_1  + e_j, e_1  + \lambda  e_j)\} \left[E(e_i^*)(1),
E(e_j^*)(\lambda)\right]^{-1}.
\end{eqnarray*}
Now by [\cite{4}, Lemma 5.2], 
\begin{eqnarray*}
E(e_i^*)(-2\lambda)^{bot}   &  =  &   S_r(e_1  -   e_j,  e_1   -  e_i)
S_r((1+\lambda)e_1 + e_j, e_1 - \lambda e_j) \\ & & S_r(e_1 - e_j, e_1
+   e_i)    S_r((1-\lambda)e_1   +    e_j,   e_1   +    \lambda   e_j)\\
& & \left[E(e_j^*)(\lambda), E(e_i^*)(1)\right]
\end{eqnarray*}
as required.
\hfill{$\Box$}

\begin{cor} \label{3.3} {\rm $($\cite{6}, Lemma 4.9, Proposition 4.10$)$}
Let $R$  be a commutative  ring with $1$  in which $2$  is invertible.
For $2 \le i \le r+1$, 
$$\mbox{the matrix of}~T_X = \left\{\begin{array}{ll} oe_{\pi(1)i}(\lambda) & \mbox{if $X =
E(e_i^*)^{bot}(-\lambda)$}\\  oe_{i\pi(1)}(-\lambda) &  \mbox{if  $X =
E(e_i)^{top}(-\lambda)$}\\   oe_{1i}(\lambda)    &   \mbox{if   $X   =
E(e_i^*)^{top}(-\lambda)$}\\   oe_{i1}(-\lambda)  &   \mbox{if   $X  =
E(e_i)^{bot}(-\lambda)$}\\
\end{array}\right.$$
\end{cor}
Proof:
By Lemma~\ref{3.1}, the matrix $A$ of $T_{S_r(e_1 - e_j, e_1 - e_i)}$ is given
by
\begin{eqnarray*}
A  &  = &  \left(  I  - \pmatrix{(e_1  -  e_j)^T  \cr  (e_1 -  e_i)^T}
\pmatrix{e_1 -  e_i & e_1  - e_j} \right) \left(I  - \pmatrix{e_1^T\cr
e_1^T}\pmatrix{e_1 & e_1}\right)\\ &  = & \pmatrix{I+e_{1i} - e_{j1} -
e_{ji} & e_{1j} -  e_{j1} - e_{jj} \cr e_{1i} - e_{i1}  - e_{ii} & I +
e_{1j} - e_{i1} - e_{ij}}.
\end{eqnarray*}
Similarly, the matrix $B$ of  $T_{S_r((1+\lambda)e_1 + e_j, e_1 - \lambda
e_j)}$ is $B = \pmatrix{B_{11} & B_{12}\cr B_{21}&B_{22}}$,
where 
\begin{eqnarray*}
B_{11} & = & I+\lambda(\lambda+2)e_{11} + \lambda(1+\lambda)e_{1j} +
(1+\lambda)e_{j1}+ \lambda e_{jj},\\
B_{12} & = & \lambda  e_{11} - (1+\lambda)e_{1j} + e_{j1} - e_{jj},\\
B_{21} & = & \lambda e_{11} + \lambda e_{1j} - \lambda(1+\lambda)e_{j1} 
- \lambda^2 e_{jj}\\
B_{22} & = & I -  e_{1j} - \lambda e_{j1} +  \lambda e_{jj},
\end{eqnarray*}
the matrix  $C$ of $T_{S_r(e_1 -  e_j, e_1 + e_i)}$ is
$$C = \pmatrix{I-e_{1i}  - e_{j1} + e_{ji} & e_{1j}  - e_{j1} - e_{jj}
  \cr -e_{1i} + e_{i1} - e_{ii} &  I + e_{1j} + e_{i1} + e_{ij}}$$ and
the matrix $D$ of $T_{S_r((1-\lambda)e_1 + e_j, e_1 + \lambda e_j)}$ is
$D = \pmatrix{D_{11} & D_{12} \cr D_{21} & D_{22}}$, where
\begin{eqnarray*}
D_{11} & = & I+\lambda(\lambda-2)e_{11} + \lambda(\lambda-1)e_{1j}+ 
(1-\lambda)e_{j1} - \lambda e_{jj},\\
D_{12} & = &  -\lambda e_{11} + (\lambda-1)e_{1j} + e_{j1} - e_{jj},\\
D_{21} & = & -\lambda e_{11}  - \lambda
  e_{1j} + \lambda(1-\lambda)e_{j1} - \lambda^2  e_{jj},\\
D_{22} & = &I - e_{1j} + \lambda e_{j1} - \lambda e_{jj}
\end{eqnarray*}
Now $AB = \pmatrix{\alpha_{11} &
  \alpha_{12} \cr \alpha_{21} & \alpha_{22}}$ where
{\small
\begin{eqnarray*}
\alpha_{11} & = & I + \lambda e_{11} + \lambda e_{1j} - \lambda e_{j1}
- \lambda e_{jj} + e_{1i} - e_{ji}\\ \alpha_{12} & = & 0\\ \alpha_{21}
&  =  & e_{1i}  -  (1+2\lambda)e_{i1} -  2\lambda  e_{ij}  - e_{ii}  -
\lambda^2 e_{11} + \lambda(1-\lambda)e_{1j} - \lambda(1+\lambda)e_{j1}
- \lambda^2 e_{jj} \\ \alpha_{22} & =  & I - e_{i1} + e_{ij} - \lambda
e_{j1} + \lambda e_{jj} - \lambda e_{11} + \lambda e_{1j}
\end{eqnarray*}
}
Also  $CD   =  \pmatrix{\beta_{11}  &  \beta_{12}   \cr  \beta_{21}  &
\beta_{22}}$ where
{\small
\begin{eqnarray*}
\beta_{11} & = & I -  \lambda e_{11} - \lambda e_{1j} + \lambda e_{j1}
+ \lambda e_{jj} - e_{1i} + e_{ji}\\ \beta_{12} & = & 0\\ \beta_{21} &
=   &   -e_{1i}   -   \lambda^2  e_{11}   -   \lambda(1+\lambda)e_{1j}
(1-2\lambda)e_{i1}     -     2\lambda     e_{ij}    -     e_{ii}     +
\lambda(1-\lambda)e_{j1}  -\lambda^2 e_{jj} \\  \beta_{22} &  = &  I -
e_{ij} - e_{1j}  + \lambda e_{j1} - \lambda e_{jj}  + \lambda e_{11} +
e_{i1}.
\end{eqnarray*}
}
Thus $$ABCD =  \pmatrix{I & 0 \cr 2\lambda e_{1i}  - 2\lambda e_{i1} -
2\lambda e_{ij}  + 2\lambda  e_{ji} & I}.$$  Also by Lemma~\ref{3.1}, the
matrix $P$ of $T_{E(e_j^*)(\lambda)}$ is given by
\begin{eqnarray*}
P   &   =   &   \left(I   -  \pmatrix{e_1^T   \cr   (e_1   +   \lambda
e_j)^T}\pmatrix{e_1   +  \lambda   e_j  &   e_1}  \right)   \left(I  -
\pmatrix{e_1^T  \cr  e_1^T}\pmatrix{e_1  &  e_1}  \right)  \\  &  =  &
\pmatrix{I - \lambda e_{1j} & 0  \cr \lambda e_{j1} - \lambda e_{1j} -
\lambda^2 e_{jj} & I + \lambda e_{j1}}.
\end{eqnarray*}
Clearly $P^{-1} = \pmatrix{I +  \lambda e_{1j} & 0 \cr -\lambda e_{j1}
  +  \lambda  e_{1j}  -  \lambda^2  e_{jj} &  I  -  \lambda  e_{j1}}$, which is the matrix of $T_{E(e_j^*)(-\lambda)}$.
Similarly,   the  matrix   $Q$  of   $T_{E(e_i^*)(1)}$  and   its  inverse $Q^{-1}$ of
$T_{E(e_i^*)(-1)}$ are
$$Q =  \pmatrix{ I - e_{1i}  & 0 \cr -e_{1i}  + e_{i1} - e_{ii}  & I +
e_{i1}},  Q^{-1} = \pmatrix{  I +  e_{1i} &  0 \cr  e_{1i} -  e_{i1} -
e_{ii} & I - e_{i1}}.$$
Thus  the matrix  $$[P, Q]  =  \pmatrix{I &  0 \cr  2\lambda e_{ij}  -
2\lambda e_{ji} &  I}.$$ Hence the product of  the matrices $ABCD$ and
$[P, Q]$ is
$$\pmatrix{I &  0 \cr  2\lambda e_{1i} -  2\lambda e_{i1}  & I} =  I +
2\lambda   e_{\pi(1)i}  -   2\lambda   e_{\pi(i)1}  =   oe_{\pi(1)i}(2
\lambda).$$  Since   $\varphi$  is  a  homomorphism,   the  matrix  of
$T_{E(e_i^*)(-2\lambda)^{bot}}$ is $oe_{\pi(1)i}(2 \lambda)$.  This proves
the    first     relation.     The    second     relation    is    its
transpose-inverse.  Similarly,  one can  prove  the  third and  fourth
relations.  \hfill{$\Box $}

\begin{cor} \label{3.4}
Let $R$  be a commutative  ring with $1$  in which $2$  is invertible.
For $2 \le i \ne j \ne \pi(i) \le r+1$, 
$$\mbox{the matrix of}~T_X = \left\{\begin{array}{ll} oe_{i\pi(1)}(\lambda)oe_{i1}(\lambda)
& \mbox{if $X = E(e_i)(\lambda)$}\\
oe_{\pi(1)i}(-\lambda)oe_{1i}(-\lambda)     &     \mbox{if    $X     =
E(e_i^*)(\lambda)$}\\             oe_{1i}(\lambda)oe_{1\pi(i)}(\lambda)
\pi_{1i}(-1)     &     \mbox{if     $X    =     E(e_{1i})(\lambda)$}\\
\pi_{1i}(-1)oe_{1i}(\lambda)oe_{1\pi(i)}(\lambda)  &   \mbox{if  $X  =
E(e_{1i}^*)(\lambda)$}\\    oe_{ij}(\lambda)   &    \mbox{if    $X   =
[E(e_j^*)(\lambda)^{top}, E(e_i)(1)^{bot}]$}\\ oe_{i\pi(j)}(\lambda) &
\mbox{if    $X    =    [E(e_j)(\lambda)^{top},    E(e_i)(1)^{bot}]$}\\
oe_{\pi(i)j}(\lambda)   &  \mbox{if  $X   =  [E(e_j^*)(\lambda)^{top},
E(e_i^*)(1)^{bot}]$}.\\
\end{array}\right.$$
$($Here  $\pi_{1i}(-1)$  denote  the matrix of $T_{S_r(e_i, e_i)}$.$)$
\end{cor}

\noindent
Proof: Follows immediately from Corollary \ref{3.3}. \hfill{$\Box $}

\begin{prop} \label{3.5}
Let   $R$  be  a   commutative  ring   with  $1$   in  which   $2$  is
invertible. Then the map $\varphi : {\rm EUm}_r(R) \to {\rm EO}_{2(r+1)}(R)$ given
by $\varphi(S_r(v,w)) = T_{S_r(v,w)}$ is surjective.
\end{prop}

\noindent
Proof: Follows from Corollary \ref{3.3}.  \hfill{$\Box $}

\section{$SUm_r(R)/EUm_r(R)$ is nilpotent}

\noindent
{\bf Notation:} Let $s$ be a non-zero divisor, $SUm_r(R, s^nR)$ denote the 
subgroup of $SUm_r(R)$ consisting of matrices which are identity modulo 
($s^n$), and $EUm_r(R, s^nR)$ denote the corresponding elementary subgroup.

\begin{lemma}\label{4.1}
Let $R$ be a commutative ring with $1$. 
Let $s$ be a non-zero divisor in Jacobson radical $J(R)$ of $R$ and 
$\beta \in SUm_r(R, s^nR)$ for $n \ge 0$.
Then the matrix of the linear transformation $T_\beta$ is in 
$SO_{2(r+1)}(R, s^nR)$.  
\end{lemma}

\noindent
Proof: Since $\beta \in {SUm}_r(R, s^nR)$, $\beta = S_r(v,w)$ where 
$v \equiv e_1 {\rm ~mod~} (s^n)$ and 
$w \equiv e_1 {\rm ~mod~} (s^n)$. Let $v = (a_0, a_1, \ldots, a_r)$ and 
$w = (b_0, b_1, \ldots, b_r)$, where $a_0$ and $b_0$ are 
$\equiv 1 {\rm ~mod~} (s^n)$,
$a_i$ and $b_i$ are $\equiv 0 {\rm ~mod~} (s^n)$.  By definition, the matrix of $T_\beta$, $[T_\beta] \in {\rm SO}_{2(r+1)}(R)$ and by
Lemma~\ref{3.1}, $[T_\beta]$
{\small
\begin{eqnarray*}
 & = &\!\! \left(I_{2(r+1)} - \pmatrix{v^T \cr w^T} \pmatrix{w & v} \right)
\left(I_{2(r+1)} - \pmatrix{e_1^T \cr e_1^T} \pmatrix{e_1 & e_1} \right)\\
& = &\!\! I_{2(r+1)} - \pmatrix{v^T \cr w^T} \pmatrix{w & v} - 
\pmatrix{e_1^T \cr e_1^T} \pmatrix{e_1 & e_1} + (a_0 + b_0)\pmatrix{v^T\cr w^T}
\pmatrix{e_1 & e_1}\\
& = &\!\!\! \pmatrix{I_{r+1} - v^Tw - e_1^Te_1 + (a_0+b_0)v^Te_1 &
 - v^Tv - e_1^Te_1 + (a_0+b_0)v^Te_1 \cr
 - w^Tw - e_1^Te_1 + (a_0+b_0)w^Te_1 & 
I_{r+1} - w^Tv - e_1^Te_1 + (a_0+b_0)w^Te_1}.
\end{eqnarray*}
}
Therefore, 
\begin{eqnarray*}
[T_\beta] {\rm ~mod~} (s^n) & = &\!\!\!  
{\small \pmatrix{I_{r+1} - e_1^Te_1 - e_1^Te_1 + 2e_1^Te_1 &  
 - e_1^Te_1 - e_1^Te_1 + 2e_1^Te_1 \cr
 - e_1^Te_1 - e_1^Te_1 + 2e_1^Te_1 & 
I_{r+1} - e_1^Te_1 - e_1^Te_1 + 2e_1^Te_1}} \\
& = &\!\! I_{2(r+1)}.
\end{eqnarray*}
Hence $[T_\beta] \in SO_{2(r+1)}(R, s^nR)$. 
\hfill $\Box$

\begin{lemma} \label{4.2}
Let $R$ be a commutative ring with $1$. In $EUm_r(R[X, Y, Z])$, 
$E(c)(Z)^{tb} E(d)(X^3Y)^{tb} E(c)(-Z)^{tb}$, 
where $c = e_i$ or $e_i^*$ and $d = e_j$ or $e_j^*$ 
is a product of elementary generators in $EUm_r(R[X, Y, Z])$ 
each of which is $\equiv$  $I_{2^r}$ modulo $(X)$.
\end{lemma}

\noindent
Proof: If necessary, the reader can consult [\cite{5}, Lemma 3.1] for details.
\hfill $\Box$

\begin{lemma}\label{4.3}
Let $R$ be a commutative ring with $1$. 
Let $s$ be a non-zero divisor in Jacobson radical $J(R)$ of $R$. 
Then we can write $E(c)(1)^{bot} E(d)(s^3x)^{top} E(c)(-1)^{bot}$, 
where $c = e_i$ or $e_i^*$, $d = e_j$ or $e_j^*$ and $x \in R$, as a 
product of elementary generators in ${EUm}_r(R)$ which are 
$\equiv$ $I_{2^r}$ modulo $(s)$.
\end{lemma}

\noindent
Proof: Put $Z = 1$, $X = s$ and $Y = x$ in Lemma \ref{4.2}.
\hfill $\Box$

\begin{lemma} \label{4.4}
Let $R$ be a commutative ring with $1$. 
Let $s$ be a non-zero divisor in Jacobson radical $J(R)$ of $R$.
If $u \equiv 1 {\rm ~mod~} (s^9)$ where $u \in R$ with $u^2 = 1$, then 
$[u] \perp [u^{-1}]$ is a product of elementary generators in ${ EUm}_r(R)$ 
each of which is  $\equiv I_{2^r}$ modulo $(s)$.
\end{lemma}

\noindent
Proof: Note that,
{\small
\begin{eqnarray*}
[u] \perp [u^{-1}] & = & \{E(e_2)(1-u^{-1})^{bot}E(e_2^*)(1-u^{-1})^{bot}\}
\{E(e_2^*)(-1)^{bot}E(e_2)(-1)^{bot}\}\\
& & \{E(e_2)(1-u)^{top}E(e_2^*)(1-u)^{top}\}\{E(e_2)(1)^{bot}E(e_2^*)(1)^{bot}\}\\
& & \{E(e_2)(1-u^{-1})^{top}E(e_2^*)(1-u^{-1})^{top}\}.
\end{eqnarray*}
}
Let $u^{-1} = u = 1+s^9 x$ for some $x \in R$. Then 
{\small
\begin{eqnarray*}
[u] \perp [u^{-1}] & = & \{E(e_2)(-s^{9}x)^{bot}E(e_2^*)(-s^{9}x)^{bot}\}
\{E(e_2^*)(-1)^{bot}E(e_2)(-1)^{bot}\}\\
& & \{E(e_2)(-s^9x)^{top}E(e_2^*)(-s^9x)^{top}\}\{E(e_2)(1)^{bot}E(e_2^*)(1)^{bot}\}\\
& & \{E(e_2)(-s^{9}x)^{top}E(e_2^*)(-s^{9}x)^{top}\}\\
& = & \{E(e_2)(-s^{9}x)^{bot}E(e_2^*)(-s^{9}x)^{bot}\} \alpha 
\{E(e_2)(-s^{9}x)^{top}E(e_2^*)(-s^{9}x)^{top}\}
\end{eqnarray*}
}
where 
\begin{eqnarray*}
\alpha  & = &  \{E(e_2^*)(-1)^{bot}E(e_2)(-1)^{bot}\}
\{E(e_2)(-s^9x)^{top}E(e_2^*)(-s^9x)^{top}\}\\
& & \{E(e_2)(1)^{bot}E(e_2^*)(1)^{bot}\}\\
& = & E(e_2^*)(-1)^{bot}\{E(e_2)(-1)^{bot}E(e_2)(-s^9x)^{top}E(e_2)(1)^{bot}\}\\
& & \{E(e_2)(-1)^{bot}E(e_2^*)(-s^9x)^{top}E(e_2)(1)^{bot}\}E(e_2^*)(1)^{bot}.
\end{eqnarray*}
By Lemma \ref{4.3}, each element in the bracket is a product of 
elementary generators in $EUm_r(R)$ which are 
$\equiv$ $I_{2^r}$ modulo $(s^3)$. 
Thus $$\alpha = E(e_2^*)(-1)^{bot} \left(\prod \alpha_i 
\prod \beta_i\right)E(e_2^*)(1)^{bot},$$
where each $\alpha_i, \beta_i \in EUm_r(R)$ with each one
$\equiv I_{2^r}$ mod $(s^3)$. Also we can write, 
\begin{eqnarray*}
\alpha &=& \prod \left(E(e_2^*)(-1)^{bot} \alpha_i 
E(e_2^*)(1)^{bot}\right) \prod \left(E(e_2^*)(-1)^{bot} \beta_i E(e_2^*)(1)^{bot}\right).
\end{eqnarray*}
Again by Lemma \ref{4.3}, each element in the product of $\alpha$ is a 
product of 
elementary generators in $EUm_r(R)$ which are $\equiv$ $I_{2^r}$ modulo $(s)$. 
Thus $[u] \perp [u^{-1}]$ is a product of 
elementary generators in $EUm_r(R)$ each of which are $\equiv$ $I_{2^r}$ 
modulo $(s)$.  \hfill $\Box$

\begin{lemma} \label{4.5}
Let $R$ be a commutative ring with $1$. Let $s$ be a 
non-zero divisor in Jacobson radical $J(R)$ of $R$ and $\beta \in 
{SUm}_r(R, s^nR)$ for $n >\!\!> 9$. Then $\beta$ can be written as a product 
of elementary generators in $EUm_r(R)$ where each is $\equiv I_{2^r}$ mod $(s)$.
\end{lemma}
Proof: By Lemma \ref{4.1}, $[T_\beta] \in {SO}_{2(r+1)}(R, s^nR)$. Thus by 
[\cite{hv}, Lemma 2.2], 
$\varphi(\beta)=[T_\beta] = \varepsilon_1 \ldots \varepsilon_k$ where each 
$\varepsilon_i \in {EO}_{2(r+1)}(R)$ which is $\equiv I_{2^r}$ mod $(s)$. For sufficiently
large $n$, we may assume that each $\varepsilon_i \equiv I_{2^r}$ mod $(s^p)$ where $n > p \ge 9$.
By Proposition \ref{3.5}, $\varepsilon_i = \varphi(\varepsilon_i')$ where each 
$\varepsilon_i' \in {EUm}_r(R, s^pR)$. Thus $\varphi(\beta) = 
\varphi(\varepsilon_1' \ldots \varepsilon_k')$. Hence 
$\beta (\varepsilon_1' \ldots \varepsilon_k')^{-1} \in \ker \varphi
= Z({\rm SUm}_r(R)) \subseteq {EUm}_r(R)$. By (\cite{6}, Corollary~3.5), 
$\beta (\varepsilon_1' \ldots \varepsilon_k')^{-1} = uI_{2^r}$ where $u$ is a 
unit with $u^2 =1$.  Since $\beta$ and $\varepsilon_i'$ are $\equiv I_{2^r}$ mod
$(s^p)$ ($n > p \ge 9$), $u \equiv 1$ mod $(s^p)$. 
Therefore, by Lemma \ref{4.4}, $\beta = u \varepsilon_1' \ldots \varepsilon_k'$
is a product of elementary generators each of which is $\equiv I_{2^r}$ mod 
$(s)$.\hfill $\Box$

\begin{lemma} \label{4.6}
Let $R$ be a commutative ring with $1$ in which $2$ is invertible, 
$s \in R$ a non-zero-divisor and $a \in R$. Then for $n >\!\!>0$
and $c = e_i$, or $e_i^*$,
$$\left[E(c)\left(\frac{a}{s}X\right)^{tb}, {SUm}_r(R, s^nR) \right] 
\subseteq {EUm}_r(R[X]).$$
More generally, given $p > 0$, for $n>\!\!>0$,
$$\left[{EUm}_r(R_s[X]), {SUm}_r(R, s^nR) \right] \subseteq {EUm}_r(R[X], 
s^pR[X])$$
\end{lemma}
Proof: Let $\alpha(X) = \left[E(c)\left(\frac{a}{s}X\right), \beta \right]$ where
$\beta \in {SUm}_r(R, s^nR)$. Then $\varphi(\beta) \in 
{SO}_{2(r+1)}(R, s^nR)$, where $\varphi :
{SUm}_r(R, s^nR) \to {SO}_{2(r+1)}(R, s^nR)$ is the canonical homomorphism. 
By Corollary \ref{3.4}, 
$$\varphi(E(c)\left(\frac{a}{s}X\right)) \in EO_{2(r+1)}(R_s[X]).$$ Thus by 
[\cite{hv}, Lemma 2.4],
\begin{eqnarray*}
\varphi(\alpha(X)) &\in& [EO_{2(r+1)}(R_s[X]), SO_{2(r+1)}(R, s^nR)] \subseteq
EO_{2(r+1)}(R[X]),
\end{eqnarray*}
and hence by Proposition \ref{3.5}, there exists 
$\varepsilon \in {EUm}_r(R[X])$
such that $\varphi(\alpha (X)) = \varphi(\varepsilon)$. This implies, $\varphi(X) \varepsilon^{-1}
\in \ker \varphi \subseteq Z({SUm}_r(R[X])) \subseteq {EUm}_r(R[X])$. Hence
$\alpha(X) \in {EUm}_r(R[X])$.
\hfill $\Box$

\begin{lemma}
Let $R$ be a commutative ring with $1$ in which $2$ is invertible, 
$s \in R$ a non-zero divisor and $a \in R$. Then for $n >\!\!>0$
and $c = e_i$, or $e_i^*$,
$$\left[E(c)\left(\frac{a}{s}\right), {SUm}_r(R, s^nR) \right] 
\subseteq {EUm}_r(R).$$
More generally, $[{EUm}_r(R_s), {SUm}_r(R, s^nR)] \subseteq {EUm}_r(R)$ for 
$n>\!\!>0$.
\end{lemma}

\noindent
Proof: Put $X = 1$ in Lemma \ref{4.6}.
\hfill $\Box$

\vskip0.1in
In (\cite{6}, Corollary 4.15) the quotient group 
$SUm_r(R)/EUm_r(R)$, $r \geq 2$ was shown to be nilpotent. This was obtained 
as a consequence of the Jose--Rao Theorem in (\cite{6}, Theorem 4.14) which 
asserts that this quotient unimodular vector group is a subgroup of the 
orthogonal quotient group $SO_{2(r+1)}(R)/EO_{2(r+1)}(R)$; which has been 
shown to be nilpotent in \cite{hv}. (Also see \cite{astep} for another 
proof.) We give a direct proof of the result following Bak's methods in 
\cite{1}.

\begin{theorem}\label{nilp}
Let $R$ be a commutative noetherian ring with $1$ in which $2$ is invertible 
and let $\dim R = d$. Then the group ${SUm}_r(R)/{EUm}_r(R)$ is nilpotent of 
class $d$ for $r \ge 2$.
\end{theorem}

\noindent
Proof: Let $G = {SUm}_r(R)/{EUm}_r(R)$. We prove that $Z^d = \{1\}$.
We prove by induction on $d = \dim R$. When $d = 0$, the ring $R$ is 
Artinian, so is semilocal. Hence ${Um}_{r+1}(R) = e_1E_{r+1}(R)$ and so any 
generator $S_r(v,w)$, $\langle v, w \rangle = 1$ is in $EUm_r(R)$.

Suppose $d>0$, Let $\alpha \in Z^d$, then $\alpha = [\beta, \gamma]$, where 
$\beta \in G$ and $\gamma \in Z^{d-1}$. Let $\beta'$ be the preimage of 
$\beta$ in $SUm_r(R)$.  

Choose a non-zero-divisor $s$ in $R$ such that $\beta_s' \in {EUm}_r(R_s)$ 
(such $s$ exists as $d > 0$). Consider $\overline G = 
\frac{{SUm}_r(R/s^nR)}{{EUm}_r(R/s^nR)}$ for some $n>\!\!>0$.  
By induction,
$\overline \gamma = \{1\}$ in $\overline G$. Since $EUm_r(R)$ is normal 
in $SUm_r(R)$, by
modifying $\gamma$ we may assume that $\gamma' \in {SUm}_r(R, s^nR)$ 
where $\gamma'$ is the preimage
of $\gamma$ in $SUm_r(R, s^nR)$. Thus by Lemma \ref{4.6}, $[\beta', \gamma'] 
\in {EUm}_r(R)$. Hence $\alpha = \{1\}$ in $G$. \hfill $\Box$

\subsection*{The relative case}

In this section we deduce the relative case of Theorem \ref{nilp} from the 
absolute case. 
We use the `Excision ring' $R \oplus I$ below instead of the usual 
non-noetherian Excision ring $\ZZ \oplus I$ as is usually done due to the 
work of van der Kallen in \cite{vdkg}. 

\begin{notn}
By $(${\cite[Proposition 5.6]{4}}$)$, the elementary generators, 
$$
E(c)(x)^{top}, E(c)(x)^{bot},$$  
where $c = e_{i} ~\mbox{or}~e_{j}^{\ast}$, and for $2\leq i\leq r+1$, and with 
$x \in R$, generate the Elementary Unimodular vector group $EUm_{r}(R)$. For 
simplicity, we shall denote these by $ge_{i}(x)$ below.  
\end{notn}

\begin{theorem} \label{relnilp}
Let $R$ be a commutative noetherian ring with $1$ in which $2$ is invertible and 
with $\dim R = d$. Let $I$ be an ideal of $R$. 
Then the group $SUm_r(R, I)/EUm_r(R, I)$ is nilpotent of class $d$ for 
$r \ge 2$.
\end{theorem}

\noindent
Proof: Let $G = {SUm}_r(R, I)/{EUm}_r(R, I)$. We prove that $Z^d = \{1\}$.
We prove by induction on $d = \dim R$. When $d = 0$, the ring $R$ is Artinian, 
so is semilocal.
Hence ${Um}_{r+1}(R, I) = e_1E_{r+1}(R, I)$ and so any generator 
$S_r(v,w)$, $\langle v, w \rangle = 1$ is in $EUm_r(R, I)$.

Suppose $d>0$, Let $\alpha \in Z^d$, then $\alpha = [\beta, \gamma]$, where $\beta \in G$ and 
$\gamma \in Z^{d-1}$. We can write $\beta = \textit{Id} + \beta{'}, \gamma = \textit{Id} + \gamma{'}$ for some 
$\beta{'}, \gamma{'} \in M_{2^{r}}(I).$ Let $\alpha = \textit{Id} + \alpha{'}$ for some $\alpha{'} \in M_{2^{r}}(I).$ 
Let $\tilde{\alpha} = (\textit{Id}, \alpha{'}) \in 
SUm_{r}(R\oplus I, 0\oplus I).$ In view of $(${\cite [Lemma 3.3]{acr}}$)$, 
\begin{eqnarray*}
\tilde{\alpha} &\in& EUm_{r}(R\oplus I) \cap SUm_{r}( R\oplus I, 0\oplus I) = EUm_{r}(R\oplus I, 0\oplus I)
\end{eqnarray*}
as $\frac{R\oplus I}{0\oplus I} 
\simeq R$ is a retract of $R\oplus I$. Thus,
$$\tilde{\alpha} = \prod_{k=1}^{m}\varepsilon_{k}ge_{i_{k}}(0,a_{k})
\varepsilon_{k}^{-1},~~~~\varepsilon_{k}\in 
EUm_{r}( R\oplus I),~ a_{k}\in I.$$
\par
Now, consider the homomorphism
\begin{eqnarray*}
f:R\oplus I\longrightarrow R\\
(r,i)\longmapsto r+i.
\end{eqnarray*}
\par
This $f$ induces a map\\
$$\tilde{f}:EUm_{r}( R\oplus I, 0\oplus I)\longrightarrow EUm_{r}(R).$$
\par
Clearly,
\begin{eqnarray*} \alpha &=& \tilde{f}(\tilde{\alpha})\\ 
&=& \prod_{k=1}^{m}\gamma_{k}ge_{i_{k}}(0 + a_{k})\gamma_{k}^{-1}\\
&=& \prod_{k=1}^{m}\gamma_{k}ge_{i_{k}}(a_{k})\gamma_{k}^{-1} \in EUm_{r}( R, I);~~\mbox{since} ~ a_{k}\in I,
\end{eqnarray*}

 $ \noindent \mbox{where},~\gamma_{k} = \tilde{f}(\varepsilon_{k})$.
 \hfill $\Box$
\section{Abelian quotients 
over polynomial extensions of a local ring}

In this section we use the Quillen--Suslin Local Global Principle, 
following the ideas of A. Bak in \cite{1},  to prove that if $R = A[X]$, with 
$A$ a local ring, then the quotient Unimodular Vector group is abelian.  
(The method is similar to the one in \cite{rrss} where we had used it to 
analyse the quotients of  the linear, symplectic, and orthogonal groups.

We begin with a few simple observations.

\vskip.1in
The following observation is well known, we record it here for future use:
\begin{lemma}
 \label{0.1}
  Let $R$ be a commutative ring and $v,w \in Um_{n}(R)$ be such that $v\cdot w^{t} = 1.$ If $v = e_{1}\sigma$ for some 
  $\sigma \in E_{n}(R)$ then there exists $\varepsilon \in E_{n}(R)$ such that $v= e_{1}\varepsilon$ and 
  $w = e_{1}(\varepsilon^{-1})^{t}.$
 \end{lemma}
 
 \noindent
 Proof: In view of $(${\cite[Corollary 2.8]{8}}$)$, $w\zeta = e_{1}(\sigma^{-1})^{t}$, where $\zeta = 
 I_{n} + v^{t}(e_{1}(\sigma^{-1})^{t} - w) \in E_{n}(R).$ We see that $v\zeta^{t} = v.$ Thus $e_{1}\sigma \zeta^{t} = 
 e_{1}\sigma = v.$ Upon taking $\varepsilon = e_{1}\sigma \zeta^{t}$, we have 
$v= e_{1}\varepsilon$ and  $w = e_{1}(\varepsilon^{-1})^{t}$. \hfill $\Box$

 \begin{cor} \label{0.2}
  Let $R$ be a local ring. For $r\geq 1$, $SUm_{r}(R) = EUm_{r}(R).$
 \end{cor}
 \noindent
Proof: Let $\alpha = S_{r}(v,w) \in SUm_{r}(R).$ Since $R$ is a local ring, therefore $v = e_{1}\sigma$ for some 
$\sigma \in E_{r+1}(R).$ Since $v\cdot w^{t} = 1$, by Lemma \ref{0.1}, there exists $\varepsilon \in E_{r+1}(R)$ such that 
$v= e_{1}\varepsilon$ and $w = e_{1}(\varepsilon^{-1})^{t}.$ Thus $\alpha = S_{r}(v,w) = S_{r}(e_{1}\varepsilon, e_{1}
(\varepsilon^{-1})^{t}) 
\in EUm_{r}(R)$.
\hfill $\Box$

\begin{lemma}\label{0.3}
 Let $R$ be a local ring and $\alpha(X),~\beta(X) \in SUm_{r}(R[X])$. Then, for $r\geq 2$, the commutator, 
 $$[\alpha(X),~\beta(X)] \in [\alpha(X)\alpha(0)^{-1},~~\beta(X)\beta(0)^{-1}] EUm_{r}(R[X]).$$
\end{lemma}
\noindent
Proof: Since $R$ is a local ring, in view of Corollary \ref{0.2}, 
$SUm_{r}(R) = EUm_{r}(R)$ for all $r\geq 1$. Thus $\alpha(0), \beta(0) 
\in EUm_{r}(R)$.
 
Let $ \eta = \alpha(X)\alpha(0)^{-1},~ \tau = \beta(X)\beta(0)^{-1}.$ Then, 
\begin{eqnarray*}
[\alpha(X),~\beta(X)] & =&  [\alpha(X)\alpha(0)^{-1}\alpha(0),~~\beta(X)\beta(0)^{-1}\beta(0)]\\
 &=& \eta\alpha(0)\tau\beta(0)(\eta\alpha(0))^{-1}(\tau\beta(0))^{-1}\\
 &= &\eta\tau\eta^{-1}\tau^{-1}(\tau\eta\tau^{-1}\alpha(0)\tau\eta^{-1}\tau^{-1})
(\tau\eta\beta(0)\alpha(0)^{-1}\eta^{-1}\tau^{-1})
(\tau\beta(0)^{-1}\tau^{-1}).
\end{eqnarray*}

By $(${\cite[Corollary 4.12]{6}}$)$, $EUm_{r}(R[X])$ is a normal subgroup of 
$SUm_{r}(R[X])$ for $r\geq 2$, hence 
\begin{eqnarray*}
(\tau\eta\tau^{-1}\alpha(0)\tau\eta^{-1}\tau^{-1})& \in &  EUm_{r}(R[X],\\ 
(\tau\eta\beta(0)\alpha(0)^{-1}\eta^{-1}\tau^{-1}) &\in& EUm_{r}(R[X],\\
(\tau\beta(0)^{-1}\tau^{-1}) &\in &  EUm_{r}(R[X]).
\end{eqnarray*}
Hence the result. \hfill $\Box$

\begin{theorem}
\label{4.8}
 Let $R$ be a local ring. Then the group $\frac{SUm_{r}(R[X])}{EUm_{r}(R[X])}$ is
 an abelian group for $r\geq 2.$
\end{theorem}
\noindent
Proof: Let $\alpha(X), \beta(X) \in SUm_{r}(R[X])$, 
we need to prove $[\alpha(X), \beta(X)] \in EUm_{r}(R[X]).$ In view of 
Lemma \ref{0.3}, we may assume that $\alpha(0) = \beta(0) = \textit{Id}$.
Define, 
$$\gamma(X,T) = [\alpha(XT), \beta(X)].$$
Then for every maximal ideal $\mathfrak{m}$ of $R[X]$,
$$\gamma(X,T)_{\mathfrak{m}} = [\alpha(XT)_{\mathfrak{m}}, \beta(X)_{\mathfrak{m}}].$$
Since $\beta(X)_{\mathfrak{m}}\in SUm_{r}(R[X]_{\mathfrak{m}}) = EUm_{r}(R[X]_{\mathfrak{m}})$, and in view of the normality of 
$EUm_{r}(R[X]_{\mathfrak{m}}[T]) \trianglelefteq  SUm_{r}(R[X]_{\mathfrak{m}}[T])$, 
for $r\geq 2$, one has 
 $\gamma(X,T)_{\mathfrak{m}}\in  EUm_{r}(R[X]_{\mathfrak{m}}[T])$ and $\gamma(X,0) = \textit{Id}$. Thus by the Local-Global 
 Principle, $(${\cite[Corollary 4.11]{6}}$)$, 
 $\gamma(X, T) \in EUm_{r}(R[X,T]),$ by putting $T = 1$, one gets,
$\gamma(X,1)  = [\alpha(X), \beta(X)] \in EUm_{r}(R[X]).$
\hfill $\Box$

\begin{theorem}
 \label{4.9}
 Let $R$ be a local ring and $I$ be an ideal of $R$. Then the group 
$\frac{SUm_{r}(R[X],I[X])}{EUm_{r}(R[X],I[X])}$ is an abelian group for 
$r\geq 2.$
\end{theorem}
\noindent
Proof: Let $\alpha, \beta \in SUm_{r}(R[X],I[X]).$ 
We can write $\alpha =
\textit{Id} + \alpha^{'}, \beta = \textit{Id} + \beta^{'}$ for some $\alpha^{'}, \beta^{'} \in  M_{2^{r}}(I[X]).$ 
Let $\sigma = [\alpha, \beta] = \textit{Id} + \sigma^{'}$ for some $\sigma^{'} \in M_{2^{r}}(I[X])$.  
Let $\tilde{\sigma} = (\textit{Id}, \sigma^{'}) \in SUm_{r}( R[X]\oplus I[X], 0\oplus I[X])$. In view of 
$(${\cite [Lemma 3.3]{acr}}$)$ and Theorem \ref{4.8},\\
 $\tilde{\sigma} \in EUm_{r}(R[X]\oplus I[X]) \cap SUm_{r}( R[X]\oplus I[X], 0\oplus I[X]) = 
EUm_{r}(R[X]\oplus I[X], 0\oplus I[X])$ as $\frac{R[X]\oplus I[X]}{0\oplus I[X]} 
\simeq R[X]$ is a retract of $R[X]\oplus I[X]$. Thus,
$$\tilde{\sigma} = \prod_{k=1}^{m}\varepsilon_{k}ge_{i_{k}}(0,a_{k})\varepsilon_{k}^{-1},~~~~\varepsilon_{k}\in 
EUm_{r}(R[X]\oplus I[X]),~ a_{k}\in I[X].$$
\par
Now, consider the homomorphism
\begin{eqnarray*}
f:R[X]\oplus I[X]\longrightarrow R[X]\\
(r,i)\longmapsto r+i.
\end{eqnarray*}
\par
This $f$ induces a map\\
$$\tilde{f}:EUm_{r}(R[X]\oplus I[X], 0\oplus I[X])\longrightarrow EUm_{r}(R[X])$$
\par
Clearly,\begin{eqnarray*} \sigma &=& \tilde{f}(\tilde{\sigma})= 
\prod_{k=1}^{m}\gamma_{k}ge_{i_{k}}(0 + a_{k})\gamma_{k}^{-1}\\
&=& \prod_{k=1}^{m}\gamma_{k}ge_{i_{k}}(a_{k})\gamma_{k}^{-1} \in E(n, R, I);~~\mbox{since} ~ a_{k}\in I,
\end{eqnarray*}

 $ \noindent \mbox{where},~\gamma_{k} = \tilde{f}(\varepsilon_{k}).$
\hfill $\Box$


\end{document}